\documentclass[12pt]{article}
\usepackage{amsthm,amssymb}

\def\silent#1\par    
{\par}

\def\dense{
\itemsep=0pt\parskip=0pt
}

\def\eqref#1{(\ref{#1})} 

\def\th{^{\rm th}}



\newcount\abranum
\newcount\abraref
\newcount\mpgw
\newcount\mpgh

\def\half{{1\over 2}}

\def\boxit#1{\medskip\vbox{\hrule\hbox{\vrule\kern3pt\vbox{\kern3pt
\vbox{\advance\hsize -6.8pt{\Rd
#1}\par}\kern3pt}\kern3pt\vrule}\hrule}\medskip}

\makeatletter
\renewcommand{\@seccntformat}[1]{\@nameuse{the#1}.\quad}
\makeatother

\def\proof{\medskip\noindent
 {\bf Proof.\enspace}}

\def\Proof#1. {\medskip\noindent
 {\bf Proof #1.}~}

\def\proofend{~\ifhmode
\unskip\nobreak\hfill\vrule height5pt width4pt depth2pt\medskip\fi
\ifmmode\eqno{\vrule height5pt width4pt depth2pt}\fi}

\def\eps{\varepsilon}

\def\text#1{\quad\mbox{\rm #1}\quad}
\def\binom#1#2{{#1\choose #2}}

\def\TL{{{\L}uczak}}
\def\PE{{Erd\H{o}s}}

\newtheorem{theorem}{Theorem}
\newtheorem{lemma}[theorem]{Lemma}

\theoremstyle{definition}

\newtheorem{remark}[theorem]{Remark}

\title{On the Multi-coloured Ramsey Numbers of Cycles}
\author{
 Tomasz \L uczak
\thanks{Faculty of Mathematics and Computer Science,
Adam Mickiewicz University, 61-614 Pozna\'n, Poland,
e-mail: {\tt tomasz@amu.edu.pl}}
\thanks{The author was partially supported by the Foundation for Polish Science.}
 \and
 Mikl\'os Simonovits
  \thanks{Alfr\'ed R\'enyi Institute of Mathematics, Hungarian Academy
    of Sciences H-1053 Budapest, Re\'altanoda u. 13-15., Hungary,
    e-mail: {\tt miki@renyi.hu}} 
  \thanks{The author was partially supported by the Hungarian National 
   Science Foundation grants OTKA T 026069, T 038210, T 0234702, and T 69062.}
  \and
 Jozef Skokan
  \thanks{Department of Mathematics, London School of Economics,
    Houghton Street, London WC2A 2AE, United Kingdom, 
    e-mail: {\tt jozef@member.ams.org}} 
    }
\date{August 21, 2010}

\begin{document}

\parskip=3pt 

\maketitle

\hyphenation{mono-chro-matic}

\pagestyle{myheadings} 
\thispagestyle{empty}

\begin{abstract}
  For a graph $L$ and an integer $k\geq 2$, $R_k(L)$ denotes the
  smallest integer $N$ for which for any edge-colouring of the complete
  graph $K_N$ by $k$ colours there exists a colour $i$ for which the 
  corresponding colour class contains $L$ as a subgraph.

 Bondy and {\PE} conjectured that, for an odd cycle $C_n$ on $n$ 
 vertices, 
 $$
  R_k(C_n) = 2^{k-1}(n-1)+1 \text{for $n>3$.}
 $$ 
 They proved the case when $k=2$ and also provided an upper bound 
 $R_k(C_n)\leq (k+2)!n$. Recently, this conjecture has been verified for
 $k=3$ if $n$ is large. In this note, we prove that for every integer 
 $k\geq 4$,
 $$
  R_k(C_n)\leq k2^kn+o(n), \text{ as $n\to\infty$.}
 $$
 When $n$ is even, Yongqi, Yuansheng, Feng, and Bingxi gave a
 construction, showing that $R_k(C_n)\geq (k-1)n-2k+4.$
 Here we prove that if $n$ is even, then
 $$
  R_k(C_n)\leq kn+o(n), \text{ as $n\to\infty$.}
 $$
\end{abstract}

\section{Introduction}

 In this note we shall consider Ramsey problems connected to
 edge-colourings of ordinary graphs with $k$ colours, for a given 
 $k\geq 2$, and try to ensure monochromatic cycles of a given length. 
 We shall use the standard notation. Given a graph $G=(V, E)$, $v(G)$ 
 denotes the number of vertices and $e(G)$ the number of edges in $G$. 
 For a subset $W$ of $V$,  $G[W]$ is the subgraph of $G$ induced by the
 vertices in $W$.

 For graphs $L_1,\dots,L_k$, the Ramsey number $R(L_1,\ldots,L_k)$ is
 the minimum integer $N$ such that for any edge-colouring of the
 complete graph $K_N$ by $k$ colours there exists a colour $i$ for 
 which the $i\th$ colour class contains $L_i$ as a~subgraph. For
 $L_1=L_2=\dots=L_k=L$, we set $R_k(L):=R(L_1,\ldots,L_k)$.

 The behaviour of the Ramsey number $R(C_n,C_m)$ has been studied by 
 several authors, for example, by Bondy and Erd\H{o}s, \cite{BonErd}, 
 Faudree and Schelp, \cite{FauSch}, Rosta, \cite{Rosta}, and it is 
 completely described and well-understood. Among others, it is known 
 that 
 \[
  R_2(C_n) 
    =
  \cases{2n-1, & \text{if $n\ge5$ is odd,}\cr
  \frac{3n}{2} -1, & \text{if $n\ge6$ is even.}\cr }
 \]
 Bondy and \PE~\cite{BonErd} conjectured that $R_k(C_n)=2^{k-1}(n-1)+1$ 
 for every odd $n>3$. The conjectured extremal colouring, giving the 
 lower bound, can be easily constructed recursively: for two colours, 
 take two disjoint sets of size $n-1$, colour all the pairs within each 
 set by colour 1, and colour all the pairs joining these two sets by 
 colour 2. For $i=3,\dots, k$, take two disjoint copies of the colouring 
 for $i-1$ colours and colour all the pairs joining these two copies by
 colour $i$.  The final $k$-colouring has $2^{k-1}(n-1)$ vertices and 
 every monochromatic component has either only $n-1$ vertices or it is
 bipartite and, therefore, does not contain odd cycles.
 
 As for the upper bound for $R_k(C_n)$, $k\geq 3$, \TL\ \cite{Lucz} 
 proved that if $n$ is odd, then $R_3(C_n) = 4n+o(n)$, as $n\to\infty$.
 Later, Kohayakawa, Simonovits, and Skokan \cite{Graco,odd-cycles} 
 showed that $R_3(C_n)=4n-3$ for all odd, sufficiently large values of
 $n$. The conjecture is still open for $k\geq 4$. Bondy and Erd\H{o}s 
 \cite{BonErd} remarked that they could prove $R_k(C_n)\leq (k+2)!n$ 
 for $n$ odd. In this note we shall give an upper bound which is 
 correct up to $O(k)$ factor.

 \begin{theorem}\label{upper}
  For every $k\geq 4$ and odd $n$, 
  $$
   R_k(C_n)\leq k2^kn+o(n), \text{ as $n\to\infty$.}
  $$
 \end{theorem}

 The Ramsey number $R_k(C_n)$ behaves rather differently for even 
 values of~$n$. From \cite{FauSch} and \cite{Rosta}, we know that 
 $R_2(C_n)=3n/2-1$ and, for large even~$n$, Benevides and Skokan 
 \cite{fabricio} proved that $R_3(C_n)=2n$. Yongqi, Yuansheng, 
 Feng, and Bingxi \cite{even} gave a construction yielding 
 $$
  R_k(C_n)\geq (k-1)n-2k+4.
 $$
 Here we prove the following.
 \begin{theorem}\label{upper:even}
  For every $k\geq 2$ and even $n$, 
  $$
   R_k(C_n)\leq kn+o(n), \text{ as $n\to\infty$.}
  $$
 \end{theorem} 

 In both results above, $k$ is fixed and $n$ is large. When $n$ is fixed and $k$ is large, we know even less. For $n$ odd and $k$ large, the original Bondy-Erd\H{o}s bound from \cite{BonErd}
 $$
  2^{k-1}(n-1)+1 \leq R_k(C_n)\leq (k+2)!n
 $$
 is still the best bound we have. There are better bounds for small values of $n$. For $n=3$, it is known that
  $$
  c_1 1073^{k/6} \leq R_k(C_3)\leq c_2 k!
 $$
 with $c_2<e$, see \cite{abb,chu,wan,xxer}. For $n=5$, Li \cite{Li1} proved that
  $$
   R_k(C_5)\leq \sqrt{18^k k!}, \text{ for $k\geq 2$.}
 $$
 For a fixed even $n$, it is known that
 $$
  R_k(C_n)\leq ck^{\frac{n}{n-2}},
 $$
 where $c$ is a constant depending only on $n$. Li and Lih \cite{Li2} showed that this bound is asymptotically correct for $n=4,6,10$.
\section{Tools}

 We shall make use of the following result of Erd\H{o}s and Gallai, 
 \cite{ErdGall59}.
 \begin{theorem}
 \label{thm:eg}
  Let $n\ge 3$. For any graph $G$ with at least 
  $(n-1)(v(G)-1)/2+1$ edges, $G$ contains a cycle of length at
  least $n$.
 \end{theorem}
 As an immediate consequence we obtain the following useful 
 decomposition lemma of Figaj and {\L}uczak~\cite{FigLucz2}.
 \begin{lemma}\label{claim:TL9}
  If no non-bipartite component of a graph $G$ contains a matching
  of at least $n/2$ edges, then there exists a partition 
  $V(G)=V^1\cup V^2\cup V^3$ of the vertices of $G$ for which
  \begin{enumerate}\dense
   \item[(A)] $G$ has no edges joining $V^1\cup V^2$ and $V^3$;
   \item[(B)] the subgraph $G[V^1\cup V^2]$ is bipartite, with
   						bipartition	$(V^1, V^2)$;
   \item[(C)] the subgraph $G[V^3]$ has at most $n(|V^3|-1)/2$
   						edges	and each component of $G[V^3]$ is
   						non-bipartite.
  \end{enumerate}
 \end{lemma}
 Notice that Lemma \ref{claim:TL9} defines a decomposition of~$V(G)$ 
 into sets $V^1$, $V^2$, and $V^3$, and we shall call $V^3$ the 
 \emph{sparse set}. For completeness, we include its proof.

 \proof
  Assume that no non-bipartite component of $G$ contains a 
  matching of at least $n/2$ edges. Let $(V^1, V^2)$ be the bipartition 
  of the union of bipartite components of $G$, and $V^3$ be the vertex 
  set of the union of non-bipartite components. Then (A) and (B) hold, 
  and (C) follows from Theorem~\ref{thm:eg}.
 \proofend

\section{Odd cycles}

 For an integer $k\geq 1$ and a positive number $c$, let $P_k(c)$ be 
the following property: 

\begin{quote}$P_k(c)$:
{\sl for every $\varepsilon>0$ there exist a $\delta>0$ and an $n_0$ 
 such that for every odd $n>n_0$ and any graph $G$ with 
 $v(G)>(1+\varepsilon)cn$ and $e(G)\geq (1-\delta)\binom{v(G)}{2}$, any 
 $k$-edge-colouring of~$G$ has a monochromatic non-bipartite component 
 with a matching of $(n+1)/2$ edges.} 
\end{quote}

 Our proof of Theorem~\ref{upper} is based on the following consequence 
 of the Regularity Lemma \cite{SzemRegu} observed by Figaj and 
 {\L}uczak; see~\cite{FigLucz2} for a more general statement.

 \begin{lemma}
  Let a real number $c>0$ and an integer $k\geq 2$ be given. If 
  $P_k(c)$ holds, then 
  $$
   R_k(C_n)\leq (c+o(1))n, \text{ as $n\to\infty$.}
  $$
 \end{lemma} 
 Hence, Theorem \ref{upper} follows from the next lemma.

 \begin{lemma}\label{lemma:upper}
  For any integer $k\geq 4$, $P_k(k2^k)$ holds.
 \end{lemma} 

 \proof
  Given an integer $k\geq 4$ and an $\varepsilon>0$, let $n$ be a 
  sufficiently large odd integer, $\delta=\varepsilon/2^{2k+4}$ and
  $N=(1+\varepsilon)k2^{k}n$. Suppose that $G$ is a graph with 
  $v(G)\geq N$ and $e(G)\geq(1-\delta)\binom{v(G)}{2}$.
  Assume to the contrary that there exists a $k$-edge colouring of $G$
  without a monochromatic matching of $(n+1)/2$ edges in a 
  non-bipartite component. We may also assume that $\eps < 1$ and
  $v(G)=N$. Indeed, if $v(G)>N$ and  
  \begin{equation}\label{eq:edges}
   e(G)\geq (1-\delta)\binom{v(G)}{2},
  \end{equation}
  then, iteratively removing ($v(G)-N$ times) a vertex of minimum
  degree, we obtain a subgraph of $G$ with $N$ vertices and at least 
  $(1-\delta)\binom{N}{2}$ edges.
   
  For every colour $i$, let $G_i$ be the spanning subgraph of $G$ 
  induced by the edges coloured by $i$. Then no $G_i$ contains a 
  matching of $(n+1)/2$ edges in a non-bipartite component, otherwise
  we would have a contradiction.
 
  We apply Lemma \ref{claim:TL9} to $G_i$ for every 
  $i\in [k]:= \{1,\dots, k\}$ and obtain a partition into $V^1_i$, 
  $V^2_i$, and the sparse set $V^3_i$. For every $i\in [k]$, set 
  $X^1_i=V_i^1$ and $X^2_i=V_i^2\cup V_i^3$. Notice that there are $2^k$ 
  sets of the form $\displaystyle\bigcap_{\ell=1}^k X_\ell^{j_\ell}$, 
  where $j_\ell\in\{1,2\}$ for every $\ell$. Since $V^1_i$, $V^2_i$ and
  $V^3_i$ form a partition of $V(G)$ for every $i$, it is clear that 
  the above $2^k$ sets are pairwise disjoint and form a 
  partition of $V(G)$.
  
  The graph $G$ has $N=(1+\varepsilon)k2^{k}n$ vertices, therefore, 
  there is a choice of $j_\ell\in\{1,2\}$, $\ell=1,2,\dots, k$, such 
  that the size of the set $\displaystyle X = 
  \bigcap_{\ell=1}^k X_\ell^{j_\ell}$ is at least 
  $N/2^k=(1+\epsilon)kn>kn$.
  
  For every $i$, if there is and edge $e$ of colour $i$ in $X$, then it 
  must be contained in $V^3_i$ (by (A) and (B)). Hence, it is contained
  in an odd component (by (C)). Since there are no monochromatic 
  matchings of $(n+1)/2$ edges in the non-bipartite components, $X$ 
  contains no cycles longer than~$n$ in colour $i$, so, by Theorem 
  \ref{thm:eg}, there are at most $n(|X|-1)/2$ edges of colour $i$ with
  both endpoints in $X$. Hence, 
  \begin{equation}\label{upperBd}
   e(G[X])\leq kn(|X|-1)/2.
  \end{equation}
  On the other hand, from \eqref{eq:edges}, we have 
  \begin{equation}\label{lower}
    e(G[X])\geq \binom{|X|}{2}- \delta\binom{N}{2}.
  \end{equation}
  Comparing \eqref{upperBd} and \eqref{lower} yields 
  $$
   |X|\le kn+ \delta\frac{ N(N-1)}{|X|-1}.
  $$
  Using assumptions $\eps<1$, $\delta=\varepsilon/2^{2k+4}$, 
  $N\le k2^{k+1}$, and $|X|>kn$, we have that 
  $$
   \delta\frac{N(N-1)}{|X|-1}
    \le 
   2\delta \frac{N^2}{|X|}
    \le 
    2\delta \frac{\big(k2^{k+1}n\big)^2}{kn} 
     \le
   \frac{\eps kn}{2}.
   $$ 
   Thus,
  $$
   (1+\varepsilon)kn\leq |X|\leq kn +\frac{\varepsilon kn}{2},
  $$ 
  which is a contradiction.
 \proofend
\begin{remark}
 The methods of Figaj and \L uczak and the proof above give a 
 result slightly stronger than Theorem \ref{upper}.
 
 {\sl
  Given a natural number $k\geq 4$ and an $\varepsilon>0$, there
  exist a $\delta>0$ and an $n_0$ with the following property. Suppose 
  that $n>n_0$ is odd, $N\geq (1+\varepsilon)k2^{k}n$, and $G$ is a 
  graph with $v(G)\geq N$ and $e(G)\geq(1-\delta)\binom{v(G)}{2}$. 
  Then in any $k$-colouring of the edges of $G$, there exists a 
  monochromatic cycle $C_n$.
  }
  
  These types of theorems are not much more difficult than the ones on 
  the colourings of the complete graphs, however, these are the forms we 
  use in our applications.
\end{remark}

\section{Even cycles}
  
  In the proof of Theorem \ref{upper:even} we shall use another case of
  the lemma of Figaj and {\L}uczak (Lemma~3 in \cite{FigLucz2}).

 \begin{lemma}\label{figaj:even}
  Let a real number $c>0$ be given. If for every $\varepsilon>0$
  there exist a $\delta>0$ and an $n_0$ such that for every even $n>n_0$
  and any graph $G$ with $v(G)>(1+\varepsilon)cn$ and 
  $e(G)\geq (1-\delta)\binom{v(G)}{2}$, any $k$-edge-colouring of $G$ 
  has a monochromatic component containing a matching of $n/2$ edges,
  then $$R_k(C_n)\leq (c+o(1))n.$$
 \end{lemma}
 Now we prove Theorem \ref{upper:even}.
 
 \proof
  For an arbitrary $0<\varepsilon < 1$, consider any $k$-colouring of a 
  graph $G$ on $N> (1+\varepsilon)nk$ vertices and with at least
  $(1-\varepsilon/3)\binom{N}{2}$ edges. One of the colours must have 
  at least 
  $\frac{1}{k}(1-\varepsilon/3)\binom{N}{2}>\half n(N-1)+1$ edges, 
  so, by Theorem \ref{thm:eg}, this colour contains a cycle of length at
  least $n+1$. This implies the existence of a matching covering $n$ 
  vertices in a monochromatic component. Hence, Lemma \ref{figaj:even} 
  implies that $R_k(C_n)\leq (k+o(1))n$.
 \proofend

\providecommand{\bysame}{\leavevmode\hbox to3em{\hrulefill}\thinspace}
\providecommand{\MR}{\relax\ifhmode\unskip\space\fi MR~}
\providecommand{\MRhref}[2]{%
 \href{http://www.ams.org/mathscinet-getitem?mr=#1}{#2}
}
\providecommand{\href}[2]{#2}


\begin{thebibliography}{10}



\bibitem{abb}
H. Abbott and D. Hanson, \emph{A problem of Schur and its generalizations}, Acta Arith {\bf 20} (1972), 175--187.

\bibitem{fabricio}
F.~Benevides and J.~Skokan, 
\emph{The 3-colored {R}amsey number of even cycles}, 
J.~Combin. Theory Ser. B \textbf{99} (2009), no.~4, 690--708.  

\bibitem{BonErd}
J.~A. Bondy and P.~Erd{\H{o}}s, 
\emph{Ramsey numbers for cycles in graphs}, 
J.~Combin. Theory Ser. B \textbf{14} (1973), 46--54.

\bibitem{chu}
F.R.K. Chung,
\emph{On triangular and cyclic Ramsey numbers with $k$ colors}, Graphs and combinatorics (Proc. Capital Conf., George Washington Univ., Washington, D.C., 1973), pp. 236--242, Lecture Notes in Math., Vol. 406, Springer, Berlin 1974. 

\bibitem{ErdGall59}
P.~Erd{\H{o}}s and T.~Gallai, \emph{On maximal paths and circuits of
 graphs}, Acta Math. Acad. Sci. Hungar. \textbf{10} (1959), 337--356
(unbound insert). 
 
\bibitem{FauSch}
R.~J. Faudree and R.~H. Schelp, \emph{All {R}amsey numbers for cycles
 in graphs}, Discrete Math. \textbf{8} (1974), 313--329.

\bibitem{FigLucz2}
A.~Figaj and T.~{\L}uczak, 
 \emph{The {R}amsey number for a triple of
 large cycles}, {\tt arXiv:0709.0048v1 [math.CO]}.

\bibitem{Graco} 
Y. Kohayakawa, M.~Simonovits, and J. Skokan, \emph{The 3-colored
 Ramsey number of odd cycles}, Proceedings of GRACO 2005, pp.~397--402
(electronic), Electron. Notes Discrete Math., 19, Elsevier, Amsterdam,
2005. 

\bibitem{odd-cycles} 
Y. Kohayakawa, M.~Simonovits, and J. Skokan, \emph{The 3-colored
 Ramsey number of odd cycles}, J.~Combin. Theory Ser. B, to appear. 

\bibitem{Li1}
Y. Li, 
\emph{The multi-color Ramsey number of an odd cycle},
J. Graph Theory {\bf 62} (2009), no.~4, 324--328. 

\bibitem{Li2}
Y. Li and K. Lih, \emph{Multi-color Ramsey numbers of even cycles},
European J. Combin. {\bf 30} (2009), no. 1, 114--118. 

\bibitem{Lucz}
T.~{\L}uczak, \emph{{$R(C\sb n,C\sb n,C\sb n)\leq(4+o(1))n$}},
J. Combin. Theory Ser. B \textbf{75} (1999), no.~2, 174--187.

\bibitem{Rosta}
V.~Rosta, \emph{On a {R}amsey-type problem of {J}. {A}. {B}ondy and
 {P}. {E}rd{\H o}s. {I}, {II}}, J.~Combin. Theory Ser. B
\textbf{15} (1973), 94--104; ibid. 15 (1973), 105--120.

\bibitem{SzemRegu}
E.~Szemer{\'e}di, \emph{Regular partitions of graphs}, Probl\`emes
combinatoires et th\'eorie des graphes (Colloq. Internat. CNRS,
Univ. Orsay, Orsay, 1976), Colloq. Internat. CNRS, vol. 260, CNRS,
Paris, 1978, pp.~399--401. 

\bibitem{even}
S.~Yongqi, Y. Yuansheng, X. Feng, and L. Bingxi,
\emph{New lower bounds on the multicolor Ramsey numbers $R\sb r(C\sb {2m})$}, 
Graphs Combin. \textbf{22} (2006), no. 2, 283--288.


\bibitem{wan} 
H. Wan, \emph{Upper bounds for Ramsey numbers $R(3,3,...,3)$ and
    Schur numbers}, J. Graph Theory \textbf{26} (1997), 119--122.

\bibitem{xxer}
X. Xu, Z. Xie, G. Exoo, and S. Radziszowski, \emph{Constructive lower bounds on classical multicolor Ramsey numbers}, Electron. J. Combin. \textbf{11} (2004), no.~1, \#R35, 24pp. (electronic).
\end{thebibliography}
\end{document}